\documentclass[12pt]{article}

\usepackage[T1]{fontenc}
\usepackage[utf8]{inputenc}
\usepackage{lmodern}
\usepackage[english]{babel}

\usepackage{mathtools,amssymb,amsthm}
\usepackage{hyperref}
\usepackage{enumitem}

\allowdisplaybreaks

\newtheorem{theorem}{Theorem}

\newtheorem{remark}{Remark}

\newtheorem{definition}{Definition}
\newtheorem{proposition}{Proposition}

\usepackage{xcolor}
\usepackage{graphicx}
\usepackage{tikz}
\usetikzlibrary{decorations.pathreplacing}
\usepackage{subcaption}
\usepackage{tabularx}
\usepackage{caption}

\usepackage{geometry}
\title{On the existence of 4-regular matchstick graphs}

\author{Mike Winkler$^1$\quad Peter Dinkelacker$^2$\quad Stefan Vogel$^3$}

\date{\small
	$^1$Fakult\"at f\"ur Mathematik, Ruhr-Universit\"at Bochum, Germany, mike.winkler@ruhr-uni-bochum.de\\[2mm]
	$^2$Togostr. 79, 13351 Berlin, Germany, peter@grity.de\\[2mm]
	$^3$Raun, Dorfstr. 7, 08648 Bad Brambach, Germany, backebackekuchen16@gmail.com\\
}

\begin{document}
	\maketitle
	\begin{abstract}
		A matchstick graph is a planar unit-distance graph. We call it \emph{4-regular} if every vertex has degree 4. While examples of 4-regular matchstick graphs with fewer than 63 vertices are known only for $n \in \{52, 54, 57, 60\}$, we prove the existence of such graphs for every integer $n \geq 63$.
	\end{abstract}

	\section{Introduction}
	
	A \emph{matchstick graph} is a planar unit-distance graph, i.e., a graph drawn in the plane with non-intersecting straight edges of unit length. A matchstick graph is called \emph{4-regular} if every vertex has degree 4, and \emph{$(2;4)$-regular} if the degree of every vertex is either 2 or 4.
	
	Finite 4-regular matchstick graphs are known to exist for all numbers of vertices $n \geq 52$, with the exceptions of $n \in \{53,\allowbreak 55,\allowbreak 56,\allowbreak 58,\allowbreak 59,\allowbreak 61,\allowbreak 62\}$. The first example, consisting of 52 vertices, was introduced by Harborth in 1986 \cite{Harborth1986}. Further examples with 54, 57, 65, 67, 73, 74, 77, and 85 vertices were constructed by the authors in 2016 and 2017 \cite{Winkler2016, Winkler2017}. Whether there exists a 4-regular matchstick graph with fewer than 52 vertices, or a non-isomorphic example with 52 vertices, remains an open problem. Approximate solutions for the missing vertex numbers are discussed in \cite{WinklerAprox, VogelMGC}.
	
	In this paper, we establish the existence of 4-regular matchstick graphs for all sufficiently large numbers of vertices.
	
	\begin{theorem}\label{thm_1}
	    For every integer $n \geq 63$, there exists a 4-regular matchstick graph with $n$ vertices.
	\end{theorem}
	
	Our proof relies on a set of eleven $(2;4)$-regular matchstick graphs with orders ranging from 5 to 49, and six 4-regular matchstick graphs with orders between 64 and 74. We also present the four known examples with fewer than 63 vertices. Unless stated otherwise, the existence of one example implies the existence of variants (e.g., via flexible transformations). For a comprehensive catalog, see \cite{WinklerCat}.
	
	Harborth \cite{Harborth2002} provided an earlier proof for $n \geq 63$, excluding the values $65, 67, 73, 74, 77,$ and $85$. The $(2;4)$-regular matchstick graphs utilized in the present article typically possess exactly two vertices of degree 2 (with one exception shown in Figure~\ref{fig_5}), making them suitable building blocks for constructing 4-regular graphs.
	
	\begin{definition}
		Let $G$ be a 4-regular matchstick graph constructed from $k$ distinct $(2;4)$-regular matchstick graphs.
	\end{definition}
	
	\begin{proposition}\label{prop_vertex_sum}
		Let $k \in \mathbb{N}$ denote the number of $(2;4)$-regular matchstick graphs used to construct a 4-regular matchstick graph $G$. The number of vertices of $G$ is the sum of the vertices of the constituent graphs minus $k$.
	\end{proposition}
	
	\begin{proposition}\label{prop_edges}
		For every 4-regular matchstick graph $G = (V, E)$, the number of edges satisfies $|E| = 2|V|$.
	\end{proposition}
	
	\begin{proposition}\label{prop_flexibility}
		If a 4-regular matchstick graph is flexible, there exist infinitely many non-congruent realizations with the same number of vertices.
	\end{proposition}
	
	\begin{proposition}\label{prop_rigidity}
		Every 4-regular matchstick graph composed of 2 or 3 rigid $(2;4)$-regular matchstick graphs is rigid.
	\end{proposition}
	
	The geometric properties, rigidity, and flexibility of the graphs presented here were verified using the computer algebra system \textsc{Matchstick Graphs Calculator} (MGC) developed by S. Vogel \cite{VogelMGC}. This software, which runs in web browsers, employs the methods described in \cite{Vogel2016}.
	
	\begin{remark}\label{rem_mgc}
		(i) The MGC provides a constructive proof for each graph. We refer to the online repository \cite{VogelMGC} as the explicit coordinate proofs are too extensive for this print. (ii) The software includes an animation feature to demonstrate flexibility (see, e.g., Figure~\ref{fig_1d}). For rigid graphs, this feature can illustrate the unique realization (e.g., Figure~\ref{fig_2b}). (iii) All figures in this article are drawn to scale; edges in the drawings have uniform length.
	\end{remark}

	\section{Small examples ($n < 63$)}
	
	\begin{theorem}\label{thm_small_examples}
			4-regular matchstick graphs exist for $n \in \{52, 54, 57, 60\}$.
		\end{theorem}
		
\begin{figure}[!ht] 
	  	\centering
	  	\begin{minipage}[t]{0.45\linewidth}
  			\centering
  			\begin{tikzpicture}
  		[y=0.4pt, x=0.4pt, yscale=-1.0, xscale=1.0]
  		\draw[line width=0.01pt]
  		(3.6290,919.5051) -- (47.2650,916.0774)
  		(47.2650,916.0774) -- (90.9009,912.6496)
  		(90.9009,912.6496) -- (66.1144,876.5737)
  		(66.1144,876.5737) -- (47.2650,916.0774)
  		(47.2650,916.0774) -- (22.4785,880.0014)
  		(22.4785,880.0014) -- (3.6290,919.5051)
  		(22.4785,880.0014) -- (66.1144,876.5737)
  		(66.1144,876.5737) -- (41.3280,840.4977)
  		(41.3280,840.4977) -- (22.4785,880.0014)
  		(90.9009,912.6496) -- (134.1311,919.5051)
  		(134.1311,919.5051) -- (118.4529,878.6389)
  		(118.4529,878.6389) -- (90.9009,912.6496)
  		(41.3280,840.4977) -- (76.3140,866.8002)
  		(76.3140,866.8002) -- (81.5996,823.3502)
  		(81.5996,823.3502) -- (41.3280,840.4977)
  		(76.3140,866.8002) -- (116.5857,849.6526)
  		(116.5857,849.6526) -- (81.5996,823.3502)
  		(81.5996,823.3502) -- (121.8713,806.2026)
  		(121.8713,806.2026) -- (116.5857,849.6526)
  		(116.5857,849.6526) -- (156.8573,832.5050)
  		(156.8573,832.5050) -- (121.8713,806.2026)
  		(121.8713,806.2026) -- (162.1429,789.0550)
  		(162.1429,789.0550) -- (156.8573,832.5050)
  		(118.4529,878.6389) -- (76.3140,866.8002)
  		(3.6290,919.5051) -- (47.2650,922.9328)
  		(47.2650,922.9328) -- (90.9009,926.3605)
  		(90.9009,926.3605) -- (66.1144,962.4364)
  		(66.1144,962.4364) -- (47.2650,922.9328)
  		(47.2650,922.9328) -- (22.4785,959.0087)
  		(22.4785,959.0087) -- (3.6290,919.5051)
  		(22.4785,959.0087) -- (66.1144,962.4364)
  		(66.1144,962.4364) -- (41.3280,998.5124)
  		(41.3280,998.5124) -- (22.4785,959.0087)
  		(90.9009,926.3605) -- (134.1311,919.5051)
  		(134.1311,919.5051) -- (118.4529,960.3712)
  		(118.4529,960.3712) -- (90.9009,926.3605)
  		(41.3280,998.5124) -- (76.3140,972.2100)
  		(76.3140,972.2100) -- (81.5996,1015.6600)
  		(81.5996,1015.6600) -- (41.3280,998.5124)
  		(76.3140,972.2100) -- (116.5857,989.3575)
  		(116.5857,989.3575) -- (81.5996,1015.6599)
  		(81.5996,1015.6600) -- (121.8713,1032.8075)
  		(121.8713,1032.8075) -- (116.5857,989.3575)
  		(116.5857,989.3575) -- (156.8573,1006.5050)
  		(156.8573,1006.5051) -- (121.8713,1032.8075)
  		(121.8713,1032.8075) -- (162.1429,1049.9551)
  		(162.1429,1049.9551) -- (156.8573,1006.5050)
  		(118.4529,960.3712) -- (76.3140,972.2099)
  		(320.6568,919.5051) -- (277.0208,916.0773)
  		(277.0208,916.0774) -- (233.3849,912.6496)
  		(233.3849,912.6496) -- (258.1713,876.5737)
  		(258.1713,876.5737) -- (277.0208,916.0774)
  		(277.0208,916.0774) -- (301.8073,880.0014)
  		(301.8073,880.0014) -- (320.6568,919.5051)
  		(301.8073,880.0014) -- (258.1713,876.5737)
  		(258.1713,876.5737) -- (282.9578,840.4977)
  		(282.9578,840.4977) -- (301.8073,880.0014)
  		(233.3849,912.6496) -- (190.1547,919.5051)
  		(190.1547,919.5051) -- (205.8328,878.6389)
  		(205.8328,878.6389) -- (233.3849,912.6496)
  		(282.9578,840.4977) -- (247.9717,866.8002)
  		(247.9717,866.8002) -- (242.6862,823.3502)
  		(242.6862,823.3502) -- (282.9578,840.4977)
  		(247.9717,866.8002) -- (207.7001,849.6526)
  		(207.7001,849.6526) -- (242.6862,823.3502)
  		(242.6862,823.3502) -- (202.4145,806.2026)
  		(202.4145,806.2026) -- (207.7001,849.6526)
  		(207.7001,849.6526) -- (167.4285,832.5050)
  		(167.4285,832.5050) -- (202.4145,806.2026)
  		(202.4145,806.2026) -- (162.1429,789.0550)
  		(162.1429,789.0550) -- (167.4285,832.5050)
  		(205.8328,878.6389) -- (247.9717,866.8002)
  		(320.6568,919.5051) -- (277.0208,922.9328)
  		(277.0208,922.9328) -- (233.3849,926.3605)
  		(233.3849,926.3605) -- (258.1713,962.4364)
  		(258.1713,962.4364) -- (277.0208,922.9328)
  		(277.0208,922.9328) -- (301.8073,959.0087)
  		(301.8073,959.0087) -- (320.6568,919.5051)
  		(301.8073,959.0087) -- (258.1713,962.4364)
  		(258.1713,962.4364) -- (282.9578,998.5124)
  		(282.9578,998.5124) -- (301.8073,959.0087)
  		(233.3849,926.3605) -- (190.1547,919.5051)
  		(190.1547,919.5051) -- (205.8328,960.3712)
  		(205.8328,960.3712) -- (233.3849,926.3605)
  		(282.9578,998.5124) -- (247.9717,972.2099)
  		(247.9717,972.2100) -- (242.6862,1015.6600)
  		(242.6862,1015.6600) -- (282.9578,998.5124)
  		(247.9717,972.2100) -- (207.7001,989.3575)
  		(207.7001,989.3575) -- (242.6862,1015.6600)
  		(242.6862,1015.6600) -- (202.4145,1032.8075)
  		(202.4145,1032.8075) -- (207.7001,989.3575)
  		(207.7001,989.3575) -- (167.4285,1006.5051)
  		(167.4285,1006.5051) -- (202.4145,1032.8075)
  		(202.4145,1032.8075) -- (162.1429,1049.9551)
  		(162.1429,1049.9551) -- (167.4285,1006.5050)
  		(205.8328,960.3712) -- (247.9717,972.2099)
  		(156.8573,832.5050) -- (162.1429,875.9551)
  		(162.1429,875.9551) -- (167.4285,832.5050)
  		(156.8573,1006.5051) -- (162.1429,963.0550)
  		(162.1429,963.0551) -- (167.4285,1006.5051)
  		(118.4529,960.3712) -- (162.1429,963.0550)
  		(162.1429,963.0551) -- (205.8328,960.3712)
  		(118.4529,878.6389) -- (162.1429,875.9551)
  		(162.1429,875.9551) -- (205.8328,878.6389);
  		\end{tikzpicture}
  			\subcaption{52 vertices}
  			\label{fig_1a}
  		\end{minipage}
  		\quad\quad
  		\begin{minipage}[t]{0.45\linewidth}
	  		\centering
  			\begin{tikzpicture}
  		[y=0.4pt, x=0.4pt, yscale=-1.0, xscale=1.0]
  		\draw[line width=0.01pt]
  		(44.5790,839.9162) -- (79.8464,865.8403)
  		(79.8464,865.8403) -- (84.6636,822.3359)
  		(84.6636,822.3359) -- (44.5790,839.9162)
  		(79.8464,865.8403) -- (119.9310,848.2600)
  		(119.9310,848.2600) -- (84.6636,822.3359)
  		(84.6636,822.3359) -- (124.7483,804.7555)
  		(124.7483,804.7555) -- (119.9310,848.2600)
  		(119.9310,848.2600) -- (160.0155,830.6796)
  		(160.0155,830.6796) -- (124.7482,804.7555)
  		(124.7483,804.7555) -- (164.8329,787.1751)
  		(164.8329,787.1751) -- (160.0156,830.6796)
  		(44.5790,839.9162) -- (25.0839,879.1053)
  		(25.0839,879.1053) -- (5.5888,918.2944)
  		(68.7702,876.3940) -- (49.2751,915.5831)
  		(49.2751,915.5831) -- (92.9614,912.8718)
  		(92.9614,912.8718) -- (68.7702,876.3940)
  		(121.3986,879.5976) -- (92.9614,912.8718)
  		(121.3986,879.5976) -- (164.8329,885.0112)
  		(121.3986,879.5976) -- (164.8329,874.1841)
  		(5.5888,918.2944) -- (49.2751,915.5831)
  		(49.2751,915.5831) -- (25.0839,879.1053)
  		(25.0839,879.1053) -- (68.7702,876.3940)
  		(68.7702,876.3940) -- (44.5790,839.9162)
  		(79.8464,865.8403) -- (121.3986,879.5976)
  		(160.0155,830.6796) -- (164.8328,874.1841)
  		(164.8329,885.0112) -- (136.3947,918.2944)
  		(92.9614,912.8718) -- (136.3946,918.2944)
  		(44.5790,996.6726) -- (79.8464,970.7485)
  		(79.8464,970.7485) -- (84.6636,1014.2530)
  		(84.6636,1014.2530) -- (44.5790,996.6726)
  		(79.8464,970.7485) -- (119.9310,988.3289)
  		(119.9310,988.3289) -- (84.6636,1014.2530)
  		(84.6636,1014.2530) -- (124.7483,1031.8333)
  		(124.7483,1031.8333) -- (119.9310,988.3289)
  		(119.9310,988.3289) -- (160.0155,1005.9092)
  		(160.0155,1005.9092) -- (124.7482,1031.8333)
  		(124.7483,1031.8333) -- (164.8329,1049.4136)
  		(164.8329,1049.4136) -- (160.0156,1005.9092)
  		(44.5790,996.6726) -- (25.0839,957.4835)
  		(25.0839,957.4835) -- (5.5888,918.2944)
  		(68.7702,960.1948) -- (49.2751,921.0057)
  		(49.2751,921.0057) -- (92.9614,923.7170)
  		(92.9614,923.7170) -- (68.7702,960.1948)
  		(121.3986,956.9912) -- (92.9614,923.7170)
  		(121.3986,956.9912) -- (164.8329,951.5777)
  		(121.3986,956.9912) -- (164.8329,962.4048)
  		(5.5888,918.2944) -- (49.2751,921.0057)
  		(49.2751,921.0057) -- (25.0839,957.4835)
  		(25.0839,957.4835) -- (68.7702,960.1948)
  		(68.7702,960.1948) -- (44.5790,996.6726)
  		(79.8464,970.7485) -- (121.3986,956.9912)
  		(160.0155,1005.9092) -- (164.8328,962.4048)
  		(164.8329,951.5777) -- (136.3947,918.2944)
  		(92.9614,923.7170) -- (136.3946,918.2944)
  		(285.0867,839.9162) -- (249.8193,865.8403)
  		(249.8193,865.8403) -- (245.0021,822.3359)
  		(245.0021,822.3359) -- (285.0867,839.9162)
  		(249.8193,865.8403) -- (209.7347,848.2600)
  		(209.7346,848.2600) -- (245.0021,822.3359)
  		(245.0021,822.3359) -- (204.9175,804.7555)
  		(204.9175,804.7555) -- (209.7346,848.2600)
  		(209.7346,848.2600) -- (169.6500,830.6796)
  		(169.6500,830.6796) -- (204.9175,804.7555)
  		(204.9175,804.7555) -- (164.8329,787.1751)
  		(164.8329,787.1751) -- (169.6501,830.6796)
  		(285.0867,839.9162) -- (304.5818,879.1053)
  		(304.5818,879.1053) -- (324.0768,918.2944)
  		(260.8954,876.3940) -- (280.3905,915.5831)
  		(280.3906,915.5831) -- (236.7043,912.8718)
  		(236.7043,912.8718) -- (260.8954,876.3940)
  		(208.2671,879.5976) -- (236.7043,912.8718)
  		(208.2671,879.5976) -- (164.8328,885.0112)
  		(208.2671,879.5976) -- (164.8328,874.1841)
  		(324.0768,918.2944) -- (280.3905,915.5831)
  		(280.3906,915.5831) -- (304.5818,879.1053)
  		(304.5818,879.1053) -- (260.8954,876.3940)
  		(260.8954,876.3940) -- (285.0866,839.9162)
  		(249.8193,865.8403) -- (208.2672,879.5976)
  		(169.6500,830.6796) -- (164.8329,874.1841)
  		(164.8329,885.0112) -- (193.2711,918.2944)
  		(236.7043,912.8718) -- (193.2711,918.2944)
  		(285.0867,996.6726) -- (249.8193,970.7485)
  		(249.8193,970.7485) -- (245.0021,1014.2530)
  		(245.0021,1014.2530) -- (285.0867,996.6726)
  		(249.8193,970.7485) -- (209.7347,988.3289)
  		(209.7346,988.3289) -- (245.0021,1014.2530)
  		(245.0021,1014.2530) -- (204.9175,1031.8333)
  		(204.9175,1031.8333) -- (209.7346,988.3289)
  		(209.7346,988.3289) -- (169.6500,1005.9092)
  		(169.6500,1005.9092) -- (204.9175,1031.8333)
  		(204.9175,1031.8333) -- (164.8329,1049.4136)
  		(164.8329,1049.4136) -- (169.6501,1005.9092)
  		(285.0867,996.6726) -- (304.5818,957.4835)
  		(304.5818,957.4835) -- (324.0768,918.2944)
  		(260.8954,960.1948) -- (280.3905,921.0057)
  		(280.3906,921.0057) -- (236.7043,923.7170)
  		(236.7043,923.7170) -- (260.8954,960.1948)
  		(208.2671,956.9912) -- (236.7043,923.7170)
  		(208.2671,956.9912) -- (164.8328,951.5777)
  		(208.2671,956.9912) -- (164.8328,962.4048)
  		(324.0768,918.2944) -- (280.3905,921.0057)
  		(280.3906,921.0057) -- (304.5818,957.4835)
  		(304.5818,957.4835) -- (260.8954,960.1948)
  		(260.8954,960.1948) -- (285.0866,996.6726)
  		(249.8193,970.7485) -- (208.2672,956.9912)
  		(169.6500,1005.9092) -- (164.8329,962.4048)
  		(164.8329,951.5777) -- (193.2711,918.2944)
  		(236.7043,923.7170) -- (193.2711,918.2944);
  		\end{tikzpicture}
  			\subcaption{54 vertices}
  			\label{fig_1b}
  		\end{minipage}
	  	\par\bigskip
	  	\begin{minipage}[t]{0.45\linewidth}
  			\centering
  			\begin{tikzpicture}
  		[y=0.4pt, x=0.4pt, yscale=-1.0, xscale=1.0]
  		\draw[line width=0.01pt]
  		(143.1401,966.9862) -- (121.2550,1004.8925)
  		(132.1975,1009.3666) -- (175.9679,1009.3666)
  		(154.0827,1047.2729) -- (132.1975,1009.3666)
  		(143.1401,966.9862) -- (132.1975,1009.3666)
  		(110.3124,1047.2729) -- (132.1975,1009.3666)
  		(154.0827,1047.2729) -- (175.9679,1009.3666)
  		(175.9679,1009.3666) -- (197.8531,1047.2729)
  		(121.2550,1004.8925) -- (110.3124,1047.2729)
  		(121.2550,1004.8924) -- (90.0237,974.2257)
  		(79.0811,1016.6061) -- (121.2550,1004.8924)
  		(79.0811,1016.6061) -- (90.0237,974.2257)
  		(90.0237,974.2257) -- (47.8499,985.9393)
  		(111.9089,936.3194) -- (90.0237,974.2257)
  		(25.9647,948.0331) -- (69.7350,948.0331)
  		(69.7350,948.0331) -- (47.8499,985.9393)
  		(111.9089,936.3194) -- (69.7350,948.0331)
  		(111.9089,936.3194) -- (143.1401,966.9862)
  		(165.0253,966.9862) -- (186.9105,1004.8924)
  		(175.9679,1009.3666) -- (132.1975,1009.3666)
  		(165.0253,966.9862) -- (175.9679,1009.3666)
  		(186.9105,1004.8924) -- (197.8531,1047.2729)
  		(186.9105,1004.8924) -- (218.1417,974.2256)
  		(229.0843,1016.6061) -- (186.9105,1004.8924)
  		(229.0843,1016.6061) -- (218.1417,974.2256)
  		(218.1417,974.2256) -- (260.3156,985.9393)
  		(196.2566,936.3194) -- (218.1417,974.2257)
  		(238.4304,948.0331) -- (260.3156,985.9393)
  		(196.2566,936.3194) -- (238.4304,948.0331)
  		(165.0253,966.9862) -- (196.2566,936.3194)
  		(143.1401,966.9862) -- (154.0827,924.6058)
  		(154.0827,924.6057) -- (165.0253,966.9862)
  		(111.9089,936.3194) -- (154.0827,924.6058)
  		(25.9647,948.0331) -- (47.8499,985.9393)
  		(154.0827,924.6057) -- (196.2566,936.3194)
  		(238.4304,948.0331) -- (282.2008,948.0331)
  		(282.2008,948.0331) -- (260.3156,985.9393)
  		(207.1992,856.0327) -- (229.0843,818.1265)
  		(238.4304,825.3659) -- (260.3156,863.2722)
  		(282.2008,825.3659) -- (238.4304,825.3659)
  		(207.1992,856.0327) -- (238.4304,825.3659)
  		(260.3156,787.4597) -- (238.4304,825.3659)
  		(282.2008,825.3659) -- (260.3156,863.2722)
  		(260.3156,863.2722) -- (304.0859,863.2722)
  		(229.0843,818.1265) -- (260.3156,787.4597)
  		(229.0843,818.1265) -- (186.9105,806.4128)
  		(218.1417,775.7460) -- (229.0843,818.1265)
  		(218.1417,775.7460) -- (186.9105,806.4128)
  		(186.9105,806.4128) -- (175.9679,764.0324)
  		(165.0253,844.3190) -- (186.9105,806.4128)
  		(165.0253,844.3190) -- (154.0827,801.9386)
  		(165.0253,844.3190) -- (207.1992,856.0327)
  		(218.1417,874.9858) -- (261.9121,874.9858)
  		(218.1417,874.9858) -- (260.3156,863.2722)
  		(261.9121,874.9858) -- (304.0859,863.2722)
  		(261.9121,874.9858) -- (250.9695,917.3663)
  		(293.1434,905.6526) -- (261.9121,874.9858)
  		(250.9695,917.3663) -- (282.2008,948.0331)
  		(207.1992,917.3663) -- (238.4304,948.0331)
  		(218.1417,874.9858) -- (207.1992,917.3663)
  		(207.1992,856.0327) -- (175.9679,886.6995)
  		(175.9679,886.6995) -- (218.1417,874.9858)
  		(165.0253,844.3190) -- (175.9679,886.6995)
  		(175.9679,886.6995) -- (207.1992,917.3663)
  		(90.0237,874.9858) -- (46.2533,874.9858)
  		(47.8498,863.2722) -- (69.7350,825.3659)
  		(25.9647,825.3659) -- (47.8498,863.2722)
  		(90.0237,874.9858) -- (47.8498,863.2722)
  		(4.0795,863.2722) -- (47.8498,863.2722)
  		(25.9647,825.3659) -- (69.7350,825.3659)
  		(69.7350,825.3659) -- (47.8498,787.4597)
  		(46.2533,874.9858) -- (4.0795,863.2722)
  		(46.2533,874.9858) -- (57.1959,917.3663)
  		(15.0221,905.6526) -- (46.2533,874.9858)
  		(15.0221,905.6526) -- (57.1959,917.3663)
  		(57.1959,917.3663) -- (25.9647,948.0331)
  		(100.9663,917.3663) -- (57.1959,917.3663)
  		(100.9663,917.3663) -- (69.7350,948.0331)
  		(100.9663,917.3663) -- (90.0237,874.9858)
  		(100.9663,856.0327) -- (79.0811,818.1265)
  		(100.9663,856.0327) -- (69.7350,825.3659)
  		(79.0811,818.1265) -- (47.8498,787.4597)
  		(79.0811,818.1265) -- (121.2549,806.4128)
  		(90.0237,775.7460) -- (79.0811,818.1265)
  		(90.0237,775.7460) -- (121.2549,806.4128)
  		(121.2549,806.4128) -- (132.1975,764.0324)
  		(143.1401,844.3190) -- (121.2549,806.4128)
  		(143.1401,844.3190) -- (154.0827,801.9386)
  		(100.9663,856.0327) -- (143.1401,844.3191)
  		(143.1401,844.3190) -- (132.1975,886.6995)
  		(132.1975,886.6995) -- (100.9663,856.0327)
  		(90.0237,874.9858) -- (132.1975,886.6995)
  		(132.1975,886.6995) -- (100.9663,917.3663)
  		(250.9695,917.3663) -- (293.1434,905.6526)
  		(110.3124,1047.2729) -- (154.0827,1047.2729)
  		(154.0827,1047.2729) -- (197.8531,1047.2729)
  		(197.8531,1047.2729) -- (229.0843,1016.6061)
  		(229.0843,1016.6061) -- (260.3156,985.9393)
  		(282.2008,948.0331) -- (293.1434,905.6526)
  		(293.1434,905.6526) -- (304.0859,863.2722)
  		(304.0859,863.2722) -- (282.2008,825.3659)
  		(282.2008,825.3659) -- (260.3156,787.4597)
  		(260.3156,787.4597) -- (218.1417,775.7460)
  		(218.1417,775.7460) -- (175.9679,764.0324)
  		(132.1975,764.0323) -- (90.0237,775.7460)
  		(90.0237,775.7460) -- (47.8498,787.4597)
  		(47.8498,787.4597) -- (25.9647,825.3659)
  		(25.9647,825.3659) -- (4.0795,863.2722)
  		(4.0795,863.2722) -- (15.0221,905.6526)
  		(15.0221,905.6526) -- (25.9647,948.0331)
  		(47.8499,985.9393) -- (79.0811,1016.6061)
  		(79.0811,1016.6061) -- (110.3124,1047.2729)
  		(175.9679,764.0324) -- (132.1975,764.0324)
  		(132.1975,764.0323) -- (154.0827,801.9386)
  		(154.0827,801.9386) -- (175.9679,764.0324)
  		(250.9695,917.3663) -- (207.1992,917.3663);
  		\end{tikzpicture}
  			\subcaption{57 vertices}
  			\label{fig_1c}
  		\end{minipage}
  		\quad\quad
  		\begin{minipage}[t]{0.45\linewidth}
	  		\centering
  			\begin{tikzpicture}
  		[y=0.4pt, x=0.4pt, yscale=-1.0, xscale=1.0]
  		\draw[line width=0.01pt]
  		(170.2005,1048.5305) -- (148.3153,1010.6242)
  		(126.4301,1048.5305) -- (148.3153,1010.6242)
  		(170.2005,1048.5305) -- (192.0857,1010.6242)
  		(192.0857,1010.6242) -- (213.9708,1048.5305)
  		(148.3153,1010.6242) -- (148.3153,966.8539)
  		(148.3153,966.8539) -- (192.0857,966.8539)
  		(192.0857,966.8539) -- (192.0857,1010.6242)
  		(251.8771,1026.6453) -- (251.8771,982.8749)
  		(289.7833,1004.7601) -- (251.8771,982.8749)
  		(251.8771,1026.6453) -- (213.9708,1004.7601)
  		(213.9708,1004.7601) -- (213.9708,1048.5305)
  		(251.8771,982.8749) -- (229.9919,944.9687)
  		(192.0857,966.8539) -- (213.9708,1004.7601)
  		(28.7325,966.8539) -- (50.6177,928.9476)
  		(6.8473,928.9476) -- (50.6177,928.9476)
  		(28.7325,966.8539) -- (72.5028,966.8539)
  		(72.5028,966.8539) -- (50.6177,1004.7601)
  		(50.6177,928.9476) -- (88.5239,907.0625)
  		(88.5239,907.0625) -- (110.4091,944.9687)
  		(110.4091,944.9687) -- (72.5028,966.8539)
  		(88.5239,1026.6453) -- (126.4301,1004.7601)
  		(126.4301,1048.5305) -- (126.4301,1004.7601)
  		(88.5239,1026.6453) -- (88.5239,982.8749)
  		(88.5239,982.8749) -- (50.6177,1004.7601)
  		(126.4301,1004.7601) -- (148.3153,966.8539)
  		(148.3153,966.8539) -- (110.4091,944.9687)
  		(110.4091,944.9687) -- (88.5239,982.8749)
  		(28.7325,803.5007) -- (72.5028,803.5007)
  		(50.6177,765.5945) -- (72.5028,803.5007)
  		(28.7325,803.5007) -- (50.6177,841.4070)
  		(50.6177,841.4070) -- (6.8473,841.4070)
  		(72.5028,803.5007) -- (110.4091,825.3859)
  		(110.4091,825.3859) -- (88.5239,863.2921)
  		(88.5239,863.2921) -- (50.6177,841.4070)
  		(6.8473,885.1773) -- (44.7535,907.0625)
  		(6.8473,928.9476) -- (44.7535,907.0625)
  		(6.8473,885.1773) -- (44.7535,863.2921)
  		(44.7535,863.2921) -- (6.8473,841.4070)
  		(44.7535,907.0625) -- (88.5239,907.0625)
  		(88.5239,907.0625) -- (88.5239,863.2921)
  		(88.5239,863.2921) -- (44.7535,863.2921)
  		(170.2005,721.8241) -- (192.0857,759.7304)
  		(213.9708,721.8241) -- (192.0857,759.7304)
  		(170.2005,721.8241) -- (148.3153,759.7304)
  		(148.3153,759.7304) -- (126.4301,721.8241)
  		(192.0857,759.7304) -- (192.0857,803.5007)
  		(192.0857,803.5007) -- (148.3153,803.5007)
  		(148.3153,803.5007) -- (148.3153,759.7304)
  		(88.5239,743.7093) -- (88.5239,787.4796)
  		(50.6177,765.5945) -- (88.5239,787.4796)
  		(88.5239,743.7093) -- (126.4301,765.5945)
  		(126.4301,765.5945) -- (126.4301,721.8241)
  		(88.5239,787.4796) -- (110.4091,825.3859)
  		(110.4091,825.3859) -- (148.3153,803.5007)
  		(148.3153,803.5007) -- (126.4301,765.5945)
  		(311.6685,803.5007) -- (289.7833,841.4069)
  		(333.5537,841.4069) -- (289.7833,841.4069)
  		(311.6685,803.5007) -- (267.8981,803.5007)
  		(267.8981,803.5007) -- (289.7833,765.5945)
  		(289.7833,841.4069) -- (251.8771,863.2921)
  		(251.8771,863.2921) -- (229.9919,825.3859)
  		(229.9919,825.3859) -- (267.8981,803.5007)
  		(251.8771,743.7093) -- (213.9708,765.5945)
  		(213.9708,721.8241) -- (213.9708,765.5945)
  		(251.8771,743.7093) -- (251.8771,787.4796)
  		(251.8771,787.4796) -- (289.7833,765.5945)
  		(213.9708,765.5945) -- (192.0857,803.5007)
  		(192.0857,803.5007) -- (229.9919,825.3859)
  		(229.9919,825.3859) -- (251.8771,787.4796)
  		(311.6685,966.8539) -- (267.8981,966.8539)
  		(289.7833,1004.7601) -- (267.8981,966.8539)
  		(311.6685,966.8539) -- (289.7833,928.9476)
  		(289.7833,928.9476) -- (333.5537,928.9476)
  		(267.8981,966.8539) -- (229.9919,944.9687)
  		(229.9919,944.9687) -- (251.8771,907.0625)
  		(251.8771,907.0625) -- (289.7833,928.9476)
  		(333.5537,885.1773) -- (295.6474,863.2921)
  		(333.5537,841.4069) -- (295.6474,863.2921)
  		(333.5537,885.1773) -- (295.6474,907.0625)
  		(295.6474,907.0625) -- (333.5537,928.9476)
  		(295.6474,863.2921) -- (251.8771,863.2921)
  		(251.8771,863.2921) -- (251.8771,907.0625)
  		(251.8771,907.0625) -- (295.6474,907.0625)
  		(192.0857,966.8539) -- (229.9919,944.9687)
  		(126.4301,1048.5305) -- (170.2005,1048.5305)
  		(170.2005,1048.5305) -- (213.9708,1048.5305)
  		(213.9708,1048.5305) -- (251.8771,1026.6453)
  		(251.8771,1026.6453) -- (289.7833,1004.7601)
  		(289.7833,1004.7601) -- (311.6685,966.8539)
  		(311.6685,966.8539) -- (333.5537,928.9476)
  		(333.5537,928.9476) -- (333.5537,885.1773)
  		(333.5537,885.1773) -- (333.5537,841.4069)
  		(333.5537,841.4069) -- (311.6685,803.5007)
  		(311.6685,803.5007) -- (289.7833,765.5945)
  		(289.7833,765.5945) -- (251.8771,743.7093)
  		(251.8771,743.7093) -- (213.9708,721.8241)
  		(213.9708,721.8241) -- (170.2005,721.8241)
  		(170.2005,721.8241) -- (126.4301,721.8241)
  		(126.4301,721.8241) -- (88.5239,743.7093)
  		(88.5239,743.7093) -- (50.6177,765.5945)
  		(50.6177,765.5945) -- (28.7325,803.5007)
  		(28.7325,803.5007) -- (6.8473,841.4070)
  		(6.8473,841.4070) -- (6.8473,885.1773)
  		(6.8473,885.1773) -- (6.8473,928.9476)
  		(6.8473,928.9476) -- (28.7325,966.8539)
  		(28.7325,966.8539) -- (50.6177,1004.7601)
  		(50.6177,1004.7601) -- (88.5239,1026.6453)
  		(88.5239,1026.6453) -- (126.4301,1048.5305)
  		(72.5028,803.5007) -- (50.6177,841.4070)
  		(44.7535,863.2921) -- (44.7535,907.0625)
  		(50.6177,928.9476) -- (72.5028,966.8539)
  		(88.5239,982.8749) -- (126.4301,1004.7601)
  		(148.3153,1010.6242) -- (192.0857,1010.6242)
  		(213.9708,1004.7601) -- (251.8771,982.8749)
  		(267.8981,966.8539) -- (289.7833,928.9476)
  		(295.6474,907.0625) -- (295.6474,863.2921)
  		(289.7833,841.4069) -- (267.8981,803.5007)
  		(251.8771,787.4796) -- (213.9708,765.5945)
  		(192.0857,759.7304) -- (148.3153,759.7304)
  		(126.4301,765.5945) -- (88.5239,787.4796);
  		\end{tikzpicture}
  			\subcaption{60 vertices}
  			\label{fig_1d}
  		\end{minipage}
	  	\caption{4-regular matchstick graphs with 52, 54, 57 and 60 vertices. The graphs a, b, c are rigid; d is flexible.}
	  	\label{fig_1}
	  \end{figure}
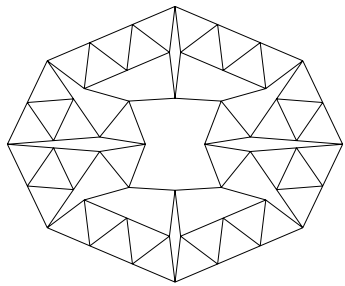
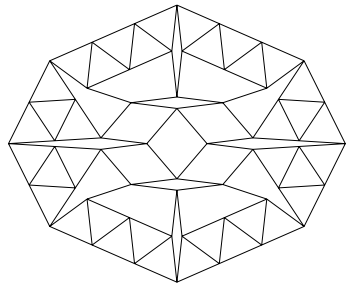
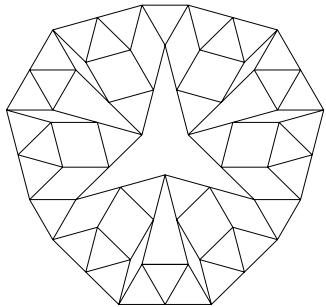
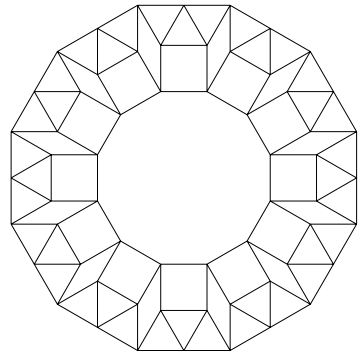

		\begin{proof}
			Figure~\ref{fig_1} depicts the known examples. Their existence is verified in \cite{VogelMGC}, with earlier proofs for the cases $n=52, 54, 57$ provided in \cite{Winkler2016, Winkler2017}.
		\end{proof}
		
	\section{$(2;4)$-regular matchstick graphs}
	
\begin{figure}[!ht]
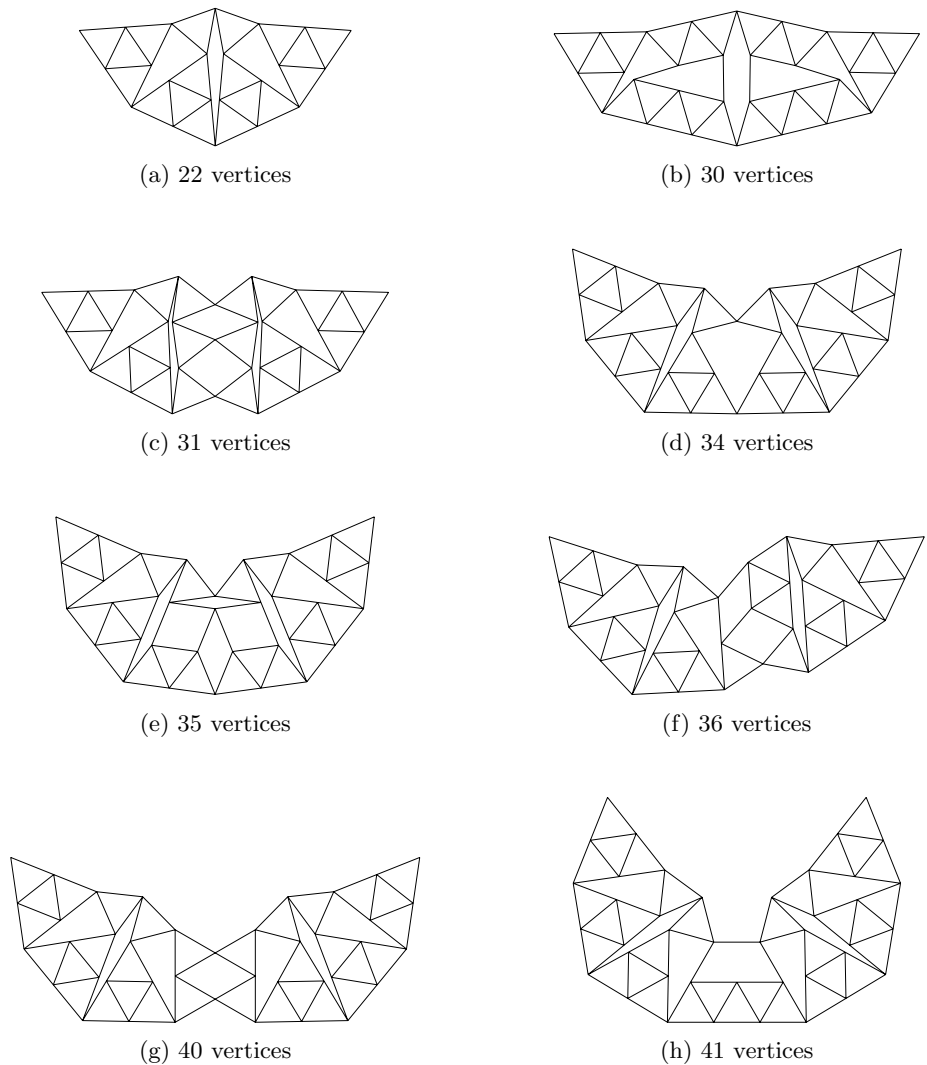
 
	  	\centering
	  	\begin{minipage}[t]{0.45\linewidth}
  			\centering

  			\subcaption{22 vertices}
  			\label{fig_2a}
  		\end{minipage}
  		\quad\quad
  		\begin{minipage}[t]{0.45\linewidth}
	  		\centering
  			%
  			\subcaption{30 vertices}
  			\label{fig_2b}
  		\end{minipage}
	  	\par\bigskip
	  	\begin{minipage}[t]{0.45\linewidth}
  			\centering
  			%
  			\subcaption{31 vertices}
  			\label{fig_2c}
  		\end{minipage}
  		\quad\quad
  		\begin{minipage}[t]{0.45\linewidth}
	  		\centering
  			%
  			\subcaption{34 vertices}
  			\label{fig_2d}
  		\end{minipage}
	  	\par\bigskip
	  	\begin{minipage}[t]{0.45\linewidth}
  			\centering
  			%
  			\subcaption{35 vertices}
  			\label{fig_2e}
  		\end{minipage}
  		\quad\quad
  		\begin{minipage}[t]{0.45\linewidth}
	  		\centering
  			%
  			\subcaption{36 vertices}
  			\label{fig_2f}
  		\end{minipage}
	  	\par\bigskip
	  	\begin{minipage}[t]{0.45\linewidth}
  			\centering
  			%
  			\subcaption{40 vertices}
  			\label{fig_2g}
  		\end{minipage}
  		\quad\quad
  		\begin{minipage}[t]{0.45\linewidth}
	  		\centering
  			%
  			\subcaption{41 vertices}
  			\label{fig_2h}
  		\end{minipage}
		\caption{$(2;4)$-regular matchstick graphs with two vertices of degree 2.}
		\label{fig_2}
	  \end{figure}

	Examples of $(2;4)$-regular matchstick graphs with fewer than 42 vertices containing exactly two vertices of degree 2 are known only for $n \in \{22,\allowbreak 30,\allowbreak 31,\allowbreak 34,\allowbreak 35,\allowbreak 36,\allowbreak 37,\allowbreak 38,\allowbreak 39,\allowbreak 40,\allowbreak 41\}$ \cite{WinklerCat}. For the construction presented in this paper, we utilize the specific examples shown in Figure~\ref{fig_2}. For each graph $G=(V,E)$ in this set, the number of edges satisfies $|E| = 2|V| - 2$.
	
	\begin{theorem}\label{thm_3}
		There exist $(2;4)$-regular matchstick graphs containing exactly two vertices of degree 2 for vertices $n \in \{22, 30, 31, 34, 35, 36, 40, 41\}$.
	\end{theorem}
	
	\begin{proof}
		Figure~\ref{fig_2} displays examples for these vertex counts. The existence of each graph is verified in \cite{VogelMGC}.
	\end{proof}	
	
	\section{Construction of 4-regular matchstick graphs with $63 \leq n \leq 120$}
	
	\begin{theorem}\label{thm_4}
		There exists a 4-regular matchstick graph for every integer $n$ with $63 \le n \le 120$, with the exceptions of $n \in \{64,\allowbreak 65,\allowbreak 67,\allowbreak 69,\allowbreak 73,\allowbreak 74\}$.
	\end{theorem}
	
	\begin{proof}
	    We consider 4-regular matchstick graphs constructed by connecting three of the eight $(2;4)$-regular graphs from Figure~\ref{fig_2}. Table~\ref{tab_1} lists all realized vertex counts $v$ and the number of distinct combinations $g$. The resulting vertex count follows from Proposition~\ref{prop_vertex_sum}.
	    
	    \begin{table}[!ht]
	    	\centering
	    	\begin{tabular}{|rr|rr|rr|rr|rr|rr|rr|rr|}
	    		\hline
	    		63&1 & 71&1 & 79&1 & 87&2 & 95&4 & 103&6 & 111&1 &119&1\\
	    		64&0 & 72&1 & 80&1 & 88&2 & 96&3 & 104&4 & 112&2 &120&1\\
	    		65&0 & 73&0 & 81&1 & 89&4 & 97&4 & 105&3 & 113&3 &&\\
	    		66&0 & 74&0 & 82&1 & 90&4 & 98&5 & 106&2 & 114&2 &&\\
	    		67&0 & 75&1 & 83&1 & 91&3 & 99&6 & 107&4 & 115&1 &&\\
	    		68&0 & 76&1 & 84&2 & 92&2 & 100&4 & 108&5 & 116&0 &&\\
	    		69&0 & 77&1 & 85&2 & 93&4 & 101&4 & 109&4 & 117&1 &&\\
	    		70&0 & 78&0 & 86&1 & 94&4 & 102&5 & 110&2 & 118&1 &&\\
	    		\hline
	    	\end{tabular}
	    	\caption{120 different possible combinations of 4-regular matchstick graphs.}
	    	\label{tab_1}
	    \end{table}
	    
	    According to Theorem~\ref{thm_3}, we can construct $8^2+\binom{8}{3}=120$ distinct rigid 4-regular matchstick graphs covering the range $63 \leq n \leq 120$ except for the values listed in Table~\ref{tab_1}.
	    The graphs 2d, 2e, 2f, 2g, and 2h can also be reflected across the axis passing through their degree-2 vertices to form rigid 4-regular matchstick graphs with 66, 68, 70, 78, and 80 vertices. Note that graph 2f requires an additional rotation by 180 degrees. The rigidity of these graphs follows from Proposition~\ref{prop_rigidity}.
	    A flexible 4-regular matchstick graph with 116 vertices is constructed by combining four copies of subgraph 2b.
	\end{proof}
	
	\begin{theorem}\label{thm_5}
		There exists a 4-regular matchstick graph for each $n \in \{64,\allowbreak 65,\allowbreak 67,\allowbreak 69,\allowbreak 73,\allowbreak 74\}$.
	\end{theorem}
	
	\begin{proof}
	    Figure~\ref{fig_3} show examples for these vertex numbers. Their existence is verified in \cite{VogelMGC}.
	\end{proof}
	
	\begin{figure}[!ht] 
	  	\centering
	  	\begin{minipage}[t]{0.45\linewidth}
  			\centering

  			\subcaption{64 vertices}
  			\label{fig_3a}
  		\end{minipage}
  		\quad\quad
  		\begin{minipage}[t]{0.45\linewidth}
	  		\centering
  			%
  			\subcaption{65 vertices}
  			\label{fig_3b}
  		\end{minipage}
  	  	 \par\bigskip
	  	\begin{minipage}[t]{0.45\linewidth}
  			\centering
  			%
  			\subcaption{67 vertices}
  			\label{fig_3c}
  		\end{minipage}
  		\quad\quad
  		\begin{minipage}[t]{0.45\linewidth}
	  		\centering
  			%
  			\subcaption{69 vertices}
  			\label{fig_3d}
  		\end{minipage}
	  	\par\bigskip
	  	\begin{minipage}[t]{0.45\linewidth}
  			\centering
  			%
  			\subcaption{73 vertices}
  			\label{fig_3e}
  		\end{minipage}
  		\quad\quad
  		\begin{minipage}[t]{0.45\linewidth}
	  		\centering
  			%
  			\subcaption{74 vertices}
  			\label{fig_3f}
  		\end{minipage}
		\caption{4-regular matchstick graphs with 64, 65, 67, 69, 73 and 74 vertices. The graphs a, d, e, f are rigid; b and c are flexible. The graph a is the only known example with 64 vertices.}
		\label{fig_3}
  	\end{figure}

	\section{Infinite families of 4-regular matchstick graphs ($n > 93$)}
	
\begin{figure}[!ht] 
	  	\centering
	  	\begin{minipage}[t]{0.29\linewidth}
  			\centering
  			\begin{tikzpicture}
  		[y=0.4pt, x=0.4pt, yscale=-1.0, xscale=1.0]
  		\draw[line width=0.01pt]
  		(4.2122,947.7114) -- (22.9573,987.2647)
  		(22.9573,987.2647) -- (41.7024,1026.8181)
  		(4.2122,947.7114) -- (47.8389,951.2544)
  		(47.8389,951.2544) -- (22.9573,987.2647)
  		(22.9573,987.2647) -- (66.5840,990.8076)
  		(66.5840,990.8076) -- (41.7024,1026.8181)
  		(47.8389,951.2544) -- (66.5840,990.8076)
  		(66.5840,990.8076) -- (91.4656,954.7973)
  		(91.4656,954.7973) -- (47.8389,951.2544)
  		(91.4656,954.7973) -- (83.9894,911.6702)
  		(48.4388,886.1359) -- (8.5502,904.1566)
  		(8.5502,904.1566) -- (44.1008,929.6908)
  		(44.1008,929.6908) -- (48.4388,886.1359)
  		(12.8881,860.6017) -- (48.4388,886.1359)
  		(48.4388,886.1359) -- (83.9894,911.6702)
  		(83.9894,911.6702) -- (44.1008,929.6908)
  		(44.1008,929.6908) -- (4.2122,947.7114)
  		(4.2122,947.7114) -- (8.5502,904.1566)
  		(91.4656,954.7973) -- (125.0767,926.7591)
  		(125.0767,926.7591) -- (83.9894,911.6702)
  		(8.5502,904.1566) -- (12.8881,860.6017)
  		(126.2749,1049.4204) -- (83.9886,1038.1192)
  		(83.9886,1038.1192) -- (41.7024,1026.8181)
  		(126.2749,1049.4204) -- (114.9189,1007.1488)
  		(114.9189,1007.1488) -- (83.9886,1038.1192)
  		(83.9886,1038.1192) -- (72.6326,995.8477)
  		(72.6326,995.8477) -- (41.7024,1026.8181)
  		(114.9189,1007.1488) -- (72.6326,995.8477)
  		(72.6326,995.8477) -- (103.5628,964.8773)
  		(103.5628,964.8773) -- (114.9189,1007.1488)
  		(103.5628,964.8773) -- (147.3311,964.4497)
  		(178.8604,994.8100) -- (168.3323,1037.2953)
  		(168.3323,1037.2953) -- (136.8030,1006.9350)
  		(136.8030,1006.9350) -- (178.8604,994.8100)
  		(210.3897,1025.1703) -- (178.8604,994.8100)
  		(178.8604,994.8100) -- (147.3311,964.4497)
  		(147.3311,964.4497) -- (136.8030,1006.9350)
  		(136.8030,1006.9350) -- (126.2749,1049.4204)
  		(126.2749,1049.4204) -- (168.3323,1037.2953)
  		(103.5628,964.8773) -- (125.0767,926.7591)
  		(125.0767,926.7591) -- (147.3311,964.4497)
  		(168.3323,1037.2953) -- (210.3897,1025.1703)
  		(83.9892,809.5330) -- (127.6159,805.9899)
  		(127.6159,805.9899) -- (171.2426,802.4469)
  		(83.9892,809.5330) -- (108.8710,845.5433)
  		(108.8710,845.5433) -- (127.6159,805.9899)
  		(127.6159,805.9899) -- (152.4977,842.0002)
  		(152.4977,842.0002) -- (171.2426,802.4469)
  		(108.8710,845.5433) -- (152.4977,842.0002)
  		(152.4977,842.0002) -- (133.7527,881.5536)
  		(133.7527,881.5536) -- (108.8710,845.5433)
  		(133.7527,881.5536) -- (92.6655,896.6427)
  		(52.7768,878.6222) -- (48.4387,835.0674)
  		(48.4387,835.0674) -- (88.3273,853.0879)
  		(88.3273,853.0879) -- (52.7768,878.6222)
  		(12.8881,860.6017) -- (52.7768,878.6222)
  		(52.7768,878.6222) -- (92.6655,896.6427)
  		(92.6655,896.6427) -- (88.3273,853.0879)
  		(88.3273,853.0879) -- (83.9892,809.5330)
  		(83.9892,809.5330) -- (48.4387,835.0674)
  		(133.7527,881.5536) -- (126.2766,924.6808)
  		(126.2766,924.6808) -- (92.6655,896.6427)
  		(48.4387,835.0674) -- (12.8881,860.6017)
  		(233.1033,864.3874) -- (202.1730,833.4171)
  		(202.1730,833.4171) -- (171.2426,802.4469)
  		(233.1033,864.3874) -- (190.8171,875.6887)
  		(190.8171,875.6887) -- (202.1730,833.4171)
  		(202.1730,833.4171) -- (159.8868,844.7185)
  		(159.8868,844.7185) -- (171.2426,802.4469)
  		(190.8171,875.6887) -- (159.8868,844.7185)
  		(159.8868,844.7185) -- (148.5309,886.9900)
  		(148.5309,886.9900) -- (190.8171,875.6887)
  		(148.5309,886.9900) -- (170.0449,925.1082)
  		(212.1024,937.2331) -- (243.6316,906.8727)
  		(243.6316,906.8727) -- (201.5741,894.7478)
  		(201.5741,894.7478) -- (212.1024,937.2331)
  		(254.1598,949.3580) -- (212.1024,937.2331)
  		(212.1024,937.2331) -- (170.0449,925.1082)
  		(170.0449,925.1082) -- (201.5741,894.7478)
  		(201.5741,894.7478) -- (233.1033,864.3874)
  		(233.1033,864.3874) -- (243.6316,906.8727)
  		(148.5309,886.9900) -- (126.2766,924.6808)
  		(126.2766,924.6808) -- (170.0449,925.1082)
  		(243.6316,906.8727) -- (254.1598,949.3580)
  		(210.3897,1025.1703) -- (232.2748,987.2640)
  		(232.2748,987.2640) -- (254.1598,949.3580)
  		(276.0452,987.2641) -- (232.2748,987.2640)
  		(276.0452,987.2641) -- (254.1601,1025.1703)
  		(254.1601,1025.1703) -- (210.3897,1025.1703)
  		(232.2748,987.2640) -- (254.1601,1025.1703)
  		(254.1601,1025.1703) -- (297.9307,1025.1703)
  		(297.9307,1025.1703) -- (276.0452,987.2641)
  		(254.1598,949.3580) -- (297.9302,949.3581)
  		(297.9302,949.3581) -- (276.0452,987.2641);
  		\end{tikzpicture}
  			\subcaption{48 vertices}
  			\label{fig_4a}
  		\end{minipage}
  		\begin{minipage}[t]{0.29\linewidth}
	  		\centering
  			\begin{tikzpicture}
  		[y=0.4pt, x=0.4pt, yscale=-1.0, xscale=1.0]
  		\draw[line width=0.01pt]
  		(3.3691,1049.7050) -- (25.2541,1011.7985)
  		(3.3686,973.8927) -- (25.2541,1011.7985)
  		(3.3691,1049.7050) -- (47.1395,1049.7047)
  		(47.1395,1049.7047) -- (25.2541,1011.7985)
  		(25.2541,1011.7985) -- (47.1392,973.8925)
  		(47.1392,973.8925) -- (3.3686,973.8927);
  		\end{tikzpicture}
  			\subcaption{5 vertices}
  			\label{fig_4b}
  		\end{minipage}
  		\begin{minipage}[t]{0.29\linewidth}
	  		\centering
  			\begin{tikzpicture}
  		[y=0.4pt, x=0.4pt, yscale=-1.0, xscale=-1.0]
  		\draw[line width=0.01pt]
  		(6.1234,935.9671) -- (21.3544,977.0020)
  		(21.3544,977.0020) -- (36.5854,1018.0369)
  		(6.1234,935.9671) -- (49.2762,943.2942)
  		(49.2762,943.2942) -- (21.3544,977.0020)
  		(21.3544,977.0020) -- (64.5071,984.3290)
  		(64.5071,984.3290) -- (36.5854,1018.0369)
  		(49.2762,943.2942) -- (64.5071,984.3290)
  		(64.5071,984.3290) -- (92.4289,950.6212)
  		(92.4289,950.6212) -- (49.2762,943.2942)
  		(92.4289,950.6212) -- (88.7351,907.0070)
  		(55.5420,878.4751) -- (14.2362,892.9552)
  		(14.2362,892.9552) -- (47.4293,921.4870)
  		(47.4293,921.4870) -- (55.5420,878.4751)
  		(22.3490,849.9433) -- (55.5420,878.4751)
  		(55.5420,878.4751) -- (88.7351,907.0070)
  		(88.7351,907.0070) -- (47.4293,921.4871)
  		(47.4293,921.4870) -- (6.1234,935.9671)
  		(6.1234,935.9671) -- (14.2362,892.9552)
  		(92.4289,950.6212) -- (128.3530,925.6151)
  		(128.3530,925.6151) -- (88.7351,907.0070)
  		(14.2362,892.9552) -- (22.3490,849.9433)
  		(118.8694,1047.9151) -- (77.7274,1032.9760)
  		(77.7274,1032.9760) -- (36.5854,1018.0369)
  		(118.8694,1047.9151) -- (111.2361,1004.8155)
  		(111.2361,1004.8155) -- (77.7274,1032.9760)
  		(77.7274,1032.9760) -- (70.0941,989.8764)
  		(70.0941,989.8764) -- (36.5854,1018.0369)
  		(111.2361,1004.8155) -- (70.0941,989.8764)
  		(70.0941,989.8764) -- (103.6027,961.7160)
  		(103.6027,961.7160) -- (111.2361,1004.8155)
  		(103.6027,961.7160) -- (147.2421,965.0999)
  		(176.0090,998.0894) -- (161.8227,1039.4970)
  		(161.8227,1039.4970) -- (133.0558,1006.5075)
  		(133.0558,1006.5075) -- (176.0090,998.0894)
  		(204.7759,1031.0789) -- (176.0090,998.0894)
  		(176.0090,998.0894) -- (147.2421,965.0999)
  		(147.2421,965.0999) -- (133.0558,1006.5075)
  		(133.0558,1006.5075) -- (118.8694,1047.9151)
  		(118.8694,1047.9151) -- (161.8227,1039.4970)
  		(103.6027,961.7160) -- (128.3530,925.6151)
  		(128.3530,925.6151) -- (147.2421,965.0999)
  		(161.8227,1039.4970) -- (204.7759,1031.0789)
  		(88.7351,792.8796) -- (131.8878,785.5525)
  		(131.8878,785.5525) -- (175.0405,778.2255)
  		(88.7351,792.8796) -- (116.6568,826.5874)
  		(116.6568,826.5874) -- (131.8878,785.5525)
  		(131.8878,785.5525) -- (159.8096,819.2604)
  		(159.8096,819.2604) -- (175.0405,778.2255)
  		(116.6568,826.5874) -- (159.8096,819.2604)
  		(159.8096,819.2604) -- (144.5786,860.2953)
  		(144.5786,860.2953) -- (116.6568,826.5874)
  		(144.5786,860.2953) -- (104.9607,878.9034)
  		(63.6548,864.4234) -- (55.5420,821.4114)
  		(55.5420,821.4114) -- (96.8479,835.8915)
  		(96.8479,835.8915) -- (63.6548,864.4233)
  		(22.3490,849.9433) -- (63.6548,864.4234)
  		(63.6548,864.4234) -- (104.9607,878.9034)
  		(104.9607,878.9034) -- (96.8479,835.8915)
  		(96.8479,835.8915) -- (88.7351,792.8796)
  		(88.7351,792.8796) -- (55.5420,821.4114)
  		(144.5786,860.2953) -- (140.8848,903.9095)
  		(140.8848,903.9095) -- (104.9607,878.9034)
  		(55.5420,821.4114) -- (22.3490,849.9433)
  		(242.0579,834.5465) -- (208.5492,806.3860)
  		(208.5492,806.3860) -- (175.0405,778.2255)
  		(242.0579,834.5465) -- (200.9159,849.4856)
  		(200.9159,849.4856) -- (208.5492,806.3860)
  		(208.5492,806.3860) -- (167.4072,821.3251)
  		(167.4072,821.3251) -- (175.0405,778.2255)
  		(200.9159,849.4856) -- (167.4072,821.3251)
  		(167.4072,821.3251) -- (159.7739,864.4247)
  		(159.7739,864.4247) -- (200.9159,849.4856)
  		(159.7739,864.4247) -- (184.5241,900.5255)
  		(227.4773,908.9436) -- (256.2442,875.9541)
  		(256.2442,875.9541) -- (213.2910,867.5360)
  		(213.2910,867.5360) -- (227.4773,908.9436)
  		(270.4306,917.3617) -- (227.4773,908.9436)
  		(227.4773,908.9436) -- (184.5241,900.5255)
  		(184.5241,900.5255) -- (213.2910,867.5360)
  		(213.2910,867.5360) -- (242.0579,834.5464)
  		(242.0579,834.5465) -- (256.2442,875.9541)
  		(159.7739,864.4247) -- (140.8848,903.9095)
  		(140.8848,903.9095) -- (184.5241,900.5255)
  		(256.2442,875.9541) -- (270.4306,917.3617)
  		(204.7759,1031.0789) -- (226.6611,993.1727)
  		(226.6611,993.1727) -- (248.5462,955.2665)
  		(248.5462,955.2665) -- (270.4306,917.3617)
  		(270.4306,917.3617) -- (292.3166,955.2675)
  		(292.3166,955.2675) -- (248.5462,955.2665)
  		(270.4314,993.1737) -- (248.5462,1031.0799)
  		(248.5462,1031.0799) -- (204.7759,1031.0789)
  		(248.5462,1031.0799) -- (226.6611,993.1727)
  		(226.6611,993.1727) -- (270.4314,993.1737)
  		(270.4314,993.1737) -- (248.5462,955.2665)
  		(248.5462,1031.0799) -- (292.3166,1031.0789)
  		(292.3166,1031.0789) -- (270.4314,993.1737);
  		\end{tikzpicture}
  			\subcaption{49 vertices}
  			\label{fig_4c}
  		\end{minipage}
	  	\caption{$(2;4)$-regular matchstick graphs with 5, 48 and 49 vertices.}
	  	\label{fig_4}
  	\end{figure}
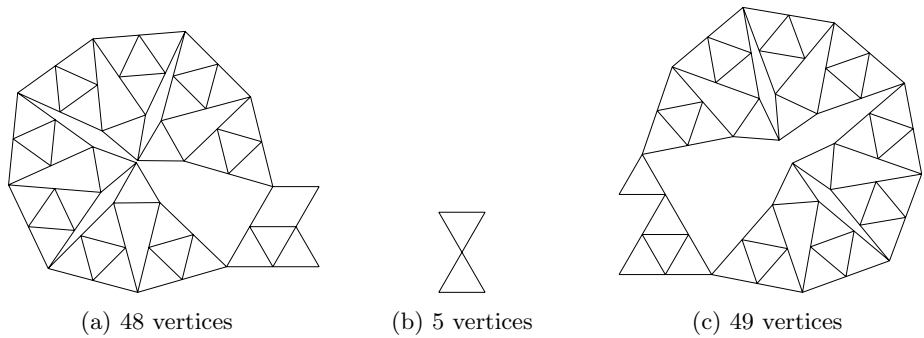

	Using graph 1a as a building block, we construct $(2;4)$-regular matchstick graphs with 48 and 49 vertices containing exactly two vertices of degree 2 (see Figure~\ref{fig_4}). The flexibility of graphs 4a and 4c allows the distance between their degree-2 vertices to be adjusted.
	
	\begin{theorem}\label{thm_6}
		There exists a 4-regular matchstick graph for each $n \in \{94, 95, 96\}$.
	\end{theorem}
	
	\begin{proof}
		The graphs 4a and 4c can be reflected across the axis through their degree-2 vertices to form 4-regular matchstick graphs with 94 and 96 vertices. Due to their flexibility, they can also be connected to each other at the degree-2 vertices to form a graph with 95 vertices. The resulting graphs are flexible.
	\end{proof}
	
	\begin{theorem}\label{thm_7}
		There exists a 4-regular matchstick graph for every integer $n \geq 97$.
	\end{theorem}
		
	\begin{proof}
		Due to its flexibility, graph 4b can be inserted into the 4-regular matchstick graphs of Theorem~\ref{thm_6} to increase the vertex count arbitrarily. Each insertion of subgraph 4b adds three vertices. This construction yields 4-regular matchstick graphs with $94+3k$, $95+3k$, and $96+3k$ vertices for any $k \in \mathbb{N}$. All resulting graphs are flexible.
	\end{proof}
	
\begin{figure}[!ht] 
		\centering 
		\begin{tikzpicture}
  	[y=0.4pt, x=0.4pt, yscale=-1.0, xscale=1.0]
  	\draw[line width=0.01pt]
  	(9.1898,945.6340) -- (27.9349,985.1873)
  	(27.9349,985.1873) -- (46.6800,1024.7406)
  	(9.1898,945.6340) -- (52.8165,949.1769)
  	(52.8165,949.1769) -- (27.9349,985.1873)
  	(27.9349,985.1873) -- (71.5616,988.7302)
  	(71.5616,988.7302) -- (46.6800,1024.7406)
  	(52.8165,949.1769) -- (71.5616,988.7302)
  	(71.5616,988.7302) -- (96.4433,952.7199)
  	(96.4433,952.7199) -- (52.8165,949.1769)
  	(96.4433,952.7199) -- (88.9670,909.5927)
  	(53.4164,884.0585) -- (13.5278,902.0791)
  	(13.5278,902.0792) -- (49.0784,927.6134)
  	(49.0784,927.6134) -- (53.4164,884.0585)
  	(17.8658,858.5243) -- (53.4164,884.0585)
  	(53.4164,884.0585) -- (88.9670,909.5927)
  	(88.9670,909.5927) -- (49.0784,927.6134)
  	(49.0784,927.6134) -- (9.1898,945.6340)
  	(9.1898,945.6340) -- (13.5278,902.0792)
  	(96.4433,952.7199) -- (130.0543,924.6817)
  	(130.0543,924.6817) -- (88.9670,909.5927)
  	(13.5278,902.0792) -- (17.8658,858.5243)
  	(131.2525,1047.3430) -- (88.9663,1036.0418)
  	(88.9663,1036.0418) -- (46.6800,1024.7406)
  	(131.2525,1047.3430) -- (119.8965,1005.0714)
  	(119.8965,1005.0714) -- (88.9663,1036.0418)
  	(88.9663,1036.0418) -- (77.6102,993.7703)
  	(77.6102,993.7703) -- (46.6800,1024.7406)
  	(119.8965,1005.0714) -- (77.6102,993.7703)
  	(77.6102,993.7703) -- (108.5405,962.7999)
  	(108.5405,962.7999) -- (119.8965,1005.0714)
  	(108.5405,962.7999) -- (152.3087,962.3723)
  	(183.8380,992.7326) -- (173.3099,1035.2179)
  	(173.3099,1035.2179) -- (141.7806,1004.8576)
  	(141.7806,1004.8576) -- (183.8380,992.7326)
  	(215.3674,1023.0929) -- (183.8380,992.7326)
  	(183.8380,992.7326) -- (152.3087,962.3723)
  	(152.3087,962.3723) -- (141.7806,1004.8576)
  	(141.7806,1004.8576) -- (131.2525,1047.3430)
  	(131.2525,1047.3430) -- (173.3099,1035.2179)
  	(108.5405,962.7999) -- (130.0543,924.6817)
  	(130.0543,924.6817) -- (152.3087,962.3723)
  	(173.3099,1035.2179) -- (215.3674,1023.0929)
  	(88.9668,807.4556) -- (132.5935,803.9125)
  	(132.5935,803.9125) -- (176.2203,800.3694)
  	(88.9668,807.4556) -- (113.8486,843.4659)
  	(113.8486,843.4659) -- (132.5935,803.9125)
  	(132.5935,803.9125) -- (157.4753,839.9228)
  	(157.4753,839.9228) -- (176.2203,800.3694)
  	(113.8486,843.4659) -- (157.4753,839.9228)
  	(157.4753,839.9228) -- (138.7303,879.4762)
  	(138.7303,879.4762) -- (113.8486,843.4659)
  	(138.7303,879.4762) -- (97.6431,894.5653)
  	(57.7544,876.5448) -- (53.4163,832.9899)
  	(53.4163,832.9899) -- (93.3050,851.0104)
  	(93.3050,851.0104) -- (57.7544,876.5448)
  	(17.8658,858.5243) -- (57.7544,876.5448)
  	(57.7544,876.5448) -- (97.6431,894.5653)
  	(97.6431,894.5653) -- (93.3050,851.0104)
  	(93.3050,851.0104) -- (88.9668,807.4556)
  	(88.9668,807.4556) -- (53.4163,832.9899)
  	(138.7303,879.4762) -- (131.2542,922.6033)
  	(131.2542,922.6033) -- (97.6431,894.5653)
  	(53.4163,832.9899) -- (17.8658,858.5243)
  	(238.0809,862.3100) -- (207.1506,831.3397)
  	(207.1506,831.3397) -- (176.2203,800.3694)
  	(238.0809,862.3100) -- (195.7947,873.6113)
  	(195.7947,873.6113) -- (207.1506,831.3397)
  	(207.1506,831.3397) -- (164.8644,842.6410)
  	(164.8644,842.6410) -- (176.2203,800.3694)
  	(195.7947,873.6113) -- (164.8644,842.6410)
  	(164.8644,842.6410) -- (153.5085,884.9126)
  	(153.5085,884.9126) -- (195.7947,873.6113)
  	(153.5085,884.9126) -- (175.0225,923.0307)
  	(217.0800,935.1556) -- (248.6092,904.7953)
  	(248.6092,904.7953) -- (206.5517,892.6703)
  	(206.5517,892.6703) -- (217.0800,935.1556)
  	(259.1374,947.2805) -- (217.0800,935.1556)
  	(217.0800,935.1556) -- (175.0225,923.0307)
  	(175.0225,923.0307) -- (206.5517,892.6703)
  	(206.5517,892.6703) -- (238.0809,862.3100)
  	(238.0809,862.3100) -- (248.6092,904.7953)
  	(153.5085,884.9126) -- (131.2542,922.6033)
  	(131.2542,922.6033) -- (175.0225,923.0307)
  	(248.6092,904.7953) -- (259.1374,947.2805)
  	(215.3674,1023.0929) -- (237.2525,985.1866)
  	(237.2525,985.1866) -- (259.1374,947.2805)
  	(281.0228,985.1867) -- (237.2525,985.1866)
  	(281.0228,985.1867) -- (259.1377,1023.0929)
  	(259.1377,1023.0929) -- (215.3674,1023.0929)
  	(237.2525,985.1866) -- (259.1377,1023.0929)
  	(259.1377,1023.0929) -- (302.9083,1023.0929)
  	(302.9083,1023.0929) -- (281.0228,985.1867)
  	(259.1374,947.2805) -- (302.9078,947.2806)
  	(302.9078,947.2806) -- (281.0228,985.1867)
  	(807.9532,927.9803) -- (792.7223,969.0152)
  	(792.7223,969.0152) -- (777.4913,1010.0500)
  	(807.9532,927.9803) -- (764.8005,935.3073)
  	(764.8005,935.3073) -- (792.7223,969.0152)
  	(792.7223,969.0152) -- (749.5695,976.3422)
  	(749.5695,976.3422) -- (777.4913,1010.0500)
  	(764.8005,935.3073) -- (749.5695,976.3422)
  	(749.5695,976.3422) -- (721.6478,942.6343)
  	(721.6478,942.6343) -- (764.8005,935.3073)
  	(721.6478,942.6343) -- (725.3416,899.0201)
  	(758.5346,870.4883) -- (799.8404,884.9684)
  	(799.8404,884.9683) -- (766.6474,913.5002)
  	(766.6474,913.5002) -- (758.5346,870.4883)
  	(791.7276,841.9564) -- (758.5346,870.4883)
  	(758.5346,870.4883) -- (725.3416,899.0201)
  	(725.3416,899.0201) -- (766.6474,913.5002)
  	(766.6474,913.5002) -- (807.9532,927.9803)
  	(807.9532,927.9803) -- (799.8404,884.9683)
  	(721.6478,942.6343) -- (685.7237,917.6283)
  	(685.7237,917.6283) -- (725.3416,899.0201)
  	(799.8404,884.9683) -- (791.7276,841.9564)
  	(695.2072,1039.9283) -- (736.3493,1024.9892)
  	(736.3493,1024.9892) -- (777.4913,1010.0500)
  	(695.2072,1039.9283) -- (702.8406,996.8287)
  	(702.8406,996.8287) -- (736.3493,1024.9892)
  	(736.3493,1024.9892) -- (743.9826,981.8896)
  	(743.9826,981.8896) -- (777.4913,1010.0500)
  	(702.8406,996.8287) -- (743.9826,981.8896)
  	(743.9826,981.8896) -- (710.4739,953.7291)
  	(710.4739,953.7291) -- (702.8406,996.8287)
  	(710.4739,953.7291) -- (666.8346,957.1130)
  	(638.0677,990.1026) -- (652.2540,1031.5102)
  	(652.2540,1031.5102) -- (681.0209,998.5207)
  	(681.0209,998.5207) -- (638.0677,990.1026)
  	(609.3008,1023.0921) -- (638.0677,990.1026)
  	(638.0677,990.1026) -- (666.8346,957.1130)
  	(666.8346,957.1130) -- (681.0209,998.5207)
  	(681.0209,998.5207) -- (695.2072,1039.9283)
  	(695.2072,1039.9283) -- (652.2540,1031.5102)
  	(710.4739,953.7291) -- (685.7237,917.6283)
  	(685.7237,917.6283) -- (666.8346,957.1130)
  	(652.2540,1031.5102) -- (609.3008,1023.0921)
  	(725.3416,784.8927) -- (682.1889,777.5657)
  	(682.1889,777.5657) -- (639.0361,770.2387)
  	(725.3416,784.8927) -- (697.4198,818.6006)
  	(697.4198,818.6006) -- (682.1889,777.5657)
  	(682.1889,777.5657) -- (654.2671,811.2735)
  	(654.2671,811.2735) -- (639.0361,770.2387)
  	(697.4198,818.6006) -- (654.2671,811.2735)
  	(654.2671,811.2735) -- (669.4981,852.3084)
  	(669.4981,852.3084) -- (697.4198,818.6006)
  	(669.4981,852.3084) -- (709.1160,870.9166)
  	(750.4218,856.4365) -- (758.5346,813.4246)
  	(758.5346,813.4246) -- (717.2288,827.9047)
  	(717.2288,827.9046) -- (750.4218,856.4365)
  	(791.7276,841.9564) -- (750.4218,856.4365)
  	(750.4218,856.4365) -- (709.1160,870.9166)
  	(709.1160,870.9166) -- (717.2288,827.9046)
  	(717.2288,827.9046) -- (725.3416,784.8927)
  	(725.3416,784.8927) -- (758.5346,813.4246)
  	(669.4981,852.3084) -- (673.1919,895.9226)
  	(673.1919,895.9226) -- (709.1160,870.9166)
  	(758.5346,813.4246) -- (791.7276,841.9564)
  	(572.0188,826.5596) -- (605.5274,798.3991)
  	(605.5274,798.3991) -- (639.0361,770.2387)
  	(572.0188,826.5596) -- (613.1608,841.4987)
  	(613.1608,841.4987) -- (605.5274,798.3991)
  	(605.5274,798.3991) -- (646.6695,813.3383)
  	(646.6695,813.3383) -- (639.0361,770.2387)
  	(613.1608,841.4987) -- (646.6695,813.3383)
  	(646.6695,813.3383) -- (654.3028,856.4379)
  	(654.3028,856.4379) -- (613.1608,841.4987)
  	(654.3028,856.4379) -- (629.5526,892.5387)
  	(586.5993,900.9568) -- (557.8324,867.9672)
  	(557.8324,867.9672) -- (600.7857,859.5491)
  	(600.7857,859.5491) -- (586.5993,900.9568)
  	(543.6461,909.3749) -- (586.5993,900.9568)
  	(586.5993,900.9568) -- (629.5526,892.5387)
  	(629.5526,892.5387) -- (600.7857,859.5491)
  	(600.7857,859.5491) -- (572.0188,826.5596)
  	(572.0188,826.5596) -- (557.8324,867.9672)
  	(654.3028,856.4379) -- (673.1919,895.9226)
  	(673.1919,895.9226) -- (629.5526,892.5387)
  	(557.8324,867.9672) -- (543.6461,909.3749)
  	(609.3008,1023.0921) -- (587.4156,985.1859)
  	(587.4156,985.1859) -- (565.5304,947.2796)
  	(565.5304,947.2796) -- (543.6461,909.3749)
  	(543.6461,909.3749) -- (521.7601,947.2806)
  	(521.7601,947.2806) -- (565.5304,947.2797)
  	(543.6453,985.1869) -- (565.5304,1023.0931)
  	(565.5304,1023.0931) -- (609.3008,1023.0931)
  	(565.5304,1023.0931) -- (587.4156,985.1859)
  	(587.4156,985.1859) -- (543.6453,985.1869)
  	(543.6453,985.1869) -- (565.5304,947.2797)
  	(565.5304,1023.0931) -- (521.7601,1023.0931)
  	(521.7601,1023.0921) -- (543.6453,985.1868);
  	
  	\draw[red,line width=0.01pt]
  	(302.9083,1023.0929) -- (324.7933,985.1863)
  	(302.9078,947.2806) -- (324.7933,985.1863)
  	(302.9083,1023.0929) -- (346.6787,1023.0926)
  	(346.6787,1023.0926) -- (324.7933,985.1863)
  	(324.7933,985.1863) -- (346.6784,947.2804)
  	(346.6784,947.2804) -- (302.9078,947.2806)
  	(477.9897,1023.0921) -- (521.7601,1023.0921)
  	(521.7601,1023.0921) -- (499.8749,985.1868)
  	(499.8749,985.1869) -- (477.9897,1023.0921)
  	(477.9897,947.2806) -- (499.8749,985.1869)
  	(499.8749,985.1869) -- (521.7601,947.2806)
  	(521.7601,947.2806) -- (477.9897,947.2806);
  	
  	\draw[line width=0.01pt]
  	(390.4485,1023.0929) -- (368.5635,985.1863)
  	(390.4490,947.2806) -- (368.5635,985.1863)
  	(390.4485,1023.0929) -- (346.6781,1023.0926)
  	(346.6781,1023.0926) -- (368.5635,985.1863)
  	(368.5635,985.1863) -- (346.6784,947.2804)
  	(346.6784,947.2804) -- (390.4490,947.2806)
  	(477.9897,1023.0921) -- (434.2194,1023.0921)
  	(434.2194,1023.0921) -- (456.1046,985.1868)
  	(456.1046,985.1869) -- (477.9897,1023.0921)
  	(477.9897,947.2806) -- (456.1046,985.1869)
  	(456.1046,985.1869) -- (434.2194,947.2806)
  	(434.2194,947.2806) -- (477.9897,947.2806);
  	
  	\path[draw=red,line width=0.8,dotted]
  	(393.4490,947.2806) -- (432.2194,947.2806)
  	(393.4485,1023.0929) -- (432.2194,1023.0921);
  	
  	\draw[black,decorate,decoration={brace,amplitude=10pt,mirror},
  	xshift=0pt,yshift=0pt] (303,1030) -- (522,1030) node [black,midway,below=15pt] {\footnotesize $n$ subgraphs 4b};
  	\end{tikzpicture}
		\caption{4-regular matchstick graph with $95+3n$ vertices.}
		\label{fig_5}
	\end{figure}
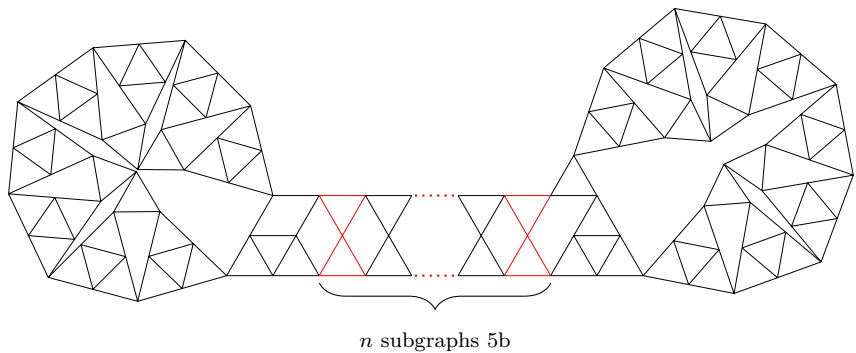

	\noindent
	Figure~\ref{fig_5} illustrates the case $95+3k$. A proof of existence is provided in \cite{VogelMGC}.
	
	\noindent
	We now present the proof of the main result.
	
	\begin{proof}[Proof of Theorem~\ref{thm_1}]
		The statement follows directly from Theorems~\ref{thm_4}, \ref{thm_5}, \ref{thm_6}, and \ref{thm_7}.
	\end{proof}

	\section{Supplementary notes}
	
	\textbf{Notes on Section 3:} Every $(2;4)$-regular matchstick graph with exactly two vertices of degree 2 is composed of smaller $(2;4)$-regular subgraphs. Consequently, multiple distinct examples may exist for the same vertex number. Figure~\ref{fig_2} displays the most symmetrical examples. With the exception of graph 2c, each example can be transformed into at least one distinct variant by reflecting and rotating the subgraphs.
	
	Table~\ref{tab_2} lists the counts of the smallest known examples ($n \leq 41$). Rotated and reflected versions of the entire graph are not counted separately. These graphs are cataloged in \cite{WinklerCat}.
	
	\begin{table}[!ht]
		\centering
		\begin{tabular}{|c|c|c|c|c|c|c|c|c|c|c|c|}
			\hline
			vertices & 22 & 30 & 31 & 34 & 35 & 36 & 37 & 38 & 39 & 40 & 41\\
			examples & 2 & 3 & 1 & 6 & 3 & 8 & 3 & 2 & 4 & 14 & 20\\
			\hline
		\end{tabular}
		\caption{Number of known examples of $(2;4)$-regular matchstick graphs with two vertices of degree 2.}
		\label{tab_2}
	\end{table}
	
	\textbf{Notes on Section 4:} The graphs in Figure~\ref{fig_3} were selected based on geometric simplicity, symmetry, and their relation to the graphs in Figure~\ref{fig_1}. Table~\ref{tab_3} provides the number of currently known 4-regular matchstick graphs for $63 \le n \le 70$. Rotated and reflected versions are not counted; flexible graphs are counted as single examples. These graphs are also presented in \cite{WinklerCat}.
	
	\begin{table}[!ht]
		\centering
		\begin{tabular}{|c|c|c|c|c|c|c|c|c|}
			\hline
			vertices & 63 & 64 & 65 & 66 & 67 & 68 & 69 & 70\\
			examples & 3 & 1 & 3 & 9 & 11 & 5 & 3 & 3\\
			\hline
		\end{tabular}
		\caption{Number of known examples of 4-regular matchstick graphs.}
		\label{tab_3}
	\end{table}

\end{document}